\documentclass[letterpaper, 10 pt, conference]{ieeeconf}
\IEEEoverridecommandlockouts
\usepackage{amsmath,amssymb,amsfonts}
\usepackage{graphicx}
\usepackage{booktabs}
\usepackage{cite}
\usepackage{capt-of}
\usepackage{xcolor}

\begin{document}

\title{Circuit-Embedded Feeder-State Reconstruction for Low-Voltage Distribution Networks under Extreme Measurement Sparsity}

\author{%
Diana Vieira Fernandes$^{1}$, Carlos Santos Silva$^{2}$, and Soummya Kar$^{3}$ 
\thanks{$^{1}$Diana Vieira Fernandes is with the Department of Engineering \& Public Policy, Carnegie Mellon University, Pittsburgh, PA, USA, and also with IN+/LARSyS, Instituto Superior T\'ecnico, Universidade de Lisboa, Lisbon, Portugal {\tt\small dianaimf@andrew.cmu.edu}}%
\thanks{$^{2}$Carlos Santos Silva is with IN+/LARSyS, Instituto Superior T\'ecnico, Universidade de Lisboa, Lisbon, Portugal {\tt\small carlos.santos.silva@tecnico.ulisboa.pt}}%
\thanks{$^{3}$Soummya Kar is with the Department of Electrical \& Computer Engineering, Carnegie Mellon University, Pittsburgh, PA, USA {\tt\small soummyak@andrew.cmu.edu}}
\thanks{This work was supported by Funda\c{c}\~ao para a Ci\^encia e a Tecnologia, Portugal, through the Carnegie Mellon Portugal Program under fellowship FCT 2025.00027.PRT and grants UIDB/50009/2025 and LA/P/0083/2020 and by ``NGS - New Generation Storage'' (C644936001-00000045), financed by PRR--Next Generation EU.}
}

\maketitle
\begin{abstract}
Low-voltage (LV) distribution feeders increasingly require internal state information, yet real-time telemetry is often limited to feeder-head active and reactive power, while downstream demand and production are described only by static contract metadata. This paper proposes a physics-based AC circuit-embedding reconstruction framework for this extreme-sparsity regime.
Feeder-head measurements, matched metadata, and class-level profiles first determine bus-level load and production priors. These priors parameterize auxiliary constant-power elements coupled to the physical feeder through impedance branches, and the feeder state is reconstructed by solving the resulting nonlinear AC circuit equilibrium. 
The resulting AC equilibrium enforces the assigned virtual-element powers, while weighted least squares (WLS) reconciles measurements and pseudo-measurements through weighted residual fitting. The method is evaluated on a real 459-bus Portuguese LV feeder model using an in-sample-calibrated synthetic benchmark of 672 fifteen-minute snapshots. 
The proposed method converged for all 672 snapshots and achieved voltage RMSE, MAE, and maximum absolute error of \(0.0033\), \(0.0014\), and \(0.0337\) p.u., respectively, outperforming both static-prior WLS baselines. Under the same constructed priors, matched-prior WLS achieved a lower RMSE of \(0.0029\,\mathrm{p.u.}\), while the proposed circuit embedding reconstructs the feeder state by solving the nonlinear AC equilibrium, bypassing WLS residual fitting entirely and requiring no specification or tuning of
pseudo-measurement variances.
\end{abstract}
%\begin{IEEEkeywords}
%Feeder-state reconstruction, circuit embedding, low-voltage networks, sparse observability, pseudo-measurements, distribution system monitoring, DER integration.
%\end{IEEEkeywords}
%=============================================================
\section{Introduction}
\label{sec:introduction}
%=============================================================
Distribution-system state estimation under limited observability has used pseudo-measurements, distributed and multi-area estimators, linearized and Kalman-filter methods, and hybrid models~\cite{dunser_ideal_2025,ferreira_kalman_2024,schubert_assessing_2025,kar_distributed_2014,kargarian_toward_2018}. Its accuracy depends strongly on prior quality and variance calibration. Backward--forward-sweep inverse power flow combines circuit calculations and ordinary least squares to infer line parameters from nodal power and voltage measurements~\cite{leal_backwardforward_2024}. Circuit-theoretic estimators embed measurements and network equations in equivalent circuits solved through modified nodal analysis~\cite{li_circuit-theoretic_2020,jovicic_enhanced_2020}. Recent generalized formulations use switch-status and PMU/SCADA measurements for joint state, topology, and bad-data estimation~\cite{li_convex_2024}. Here, topology is known, downstream telemetry is unavailable, and feeder-head active and reactive power are the only time-varying measurements. These boundary powers leave the internal injection pattern underdetermined. The proposed method first constructs bus-level load and production priors from contract metadata, profiles, and feeder-head powers, then embeds them as constant-power elements in the AC circuit equations, implemented through an augmented network with coupling branches. Its equilibrium yields physical-bus voltages, angles, and flows. The relevant comparator is pseudo-measurement weighted least squares (WLS), which encodes the same information through weighted residuals.

The main contributions are:
\begin{enumerate}
\item An AC circuit-embedding reconstruction formulation for boundary-only telemetry in which metadata-derived load and production priors are imposed as constant-power elements and the feeder state is obtained from the nonlinear AC equilibrium, requiring no pseudo-measurement variances; an augmented implementation with coupling branches provides prior-realization diagnostics and approaches direct embedding as the coupling impedances vanish, provided the coupling currents remain bounded.
\item A feeder-head-informed construction of the virtual-element powers from measured active and reactive boundary power, contracted consumption and production capacities, and class-level temporal profiles.
\item Evaluation on a real 459-bus feeder model using an in-sample-calibrated synthetic benchmark, including feeder-head consistency, two independently specified static-prior WLS baselines, direct-prior and matched-prior mechanism checks, and sensitivity to the calibration parameters \(\eta\), \(\gamma\), and \(\rho\).
\end{enumerate}
%=============================================================
\section{Problem Formulation}
\label{sec:problem}
%=============================================================
Consider an LV distribution feeder observed over discrete time indices $t \in \mathcal{T}$. The available information at each time step consists of:

\begin{enumerate}
\item \textbf{Feeder-head measurements:} active and reactive power derived from four-quadrant metering at 15-minute resolution. Specifically, $P_{\mathrm f}^{\mathrm{meas}}(t)=P_{\mathrm{import}}(t)-P_{\mathrm{export}}(t)$ and $Q_{\mathrm f}^{\mathrm{meas}}(t)=Q_1(t)+Q_2(t)-Q_3(t)-Q_4(t)$, where \(P_{\mathrm{import}}\) and \(P_{\mathrm{export}}\) are the active-power import and export channels, and \(Q_1,\dots,Q_4\) are the four-quadrant reactive-power channels, all expressed in consistent power units. Positive \(P_{\mathrm f}^{\mathrm{meas}}\) denotes net import at the feeder head, and the signs follow the meter's register convention so that positive \(Q_{\mathrm f}^{\mathrm{meas}}\) is net inductive reactive demand.
\item \textbf{Network model:} feeder topology, line parameters, transformer data, and bus connectivity.
\item \textbf{Matched contract metadata:} contracted consumption capacities and contracted or maximum production capacities associated with downstream load and production locations matched to the feeder model.
\item \textbf{Boundary-voltage assumption:} because no time-varying substation-voltage measurements are available, the feeder-boundary voltage magnitude is imposed as a fixed reference value taken from the feeder model.
\end{enumerate}
Let \(\mathcal N=\{0,1,\ldots,n\}\) denote the set of physical feeder buses, where bus \(0\) is the physical substation bus. The feeder state is represented as
\begin{equation}
\mathbf{x}(t)
=
\bigl[\, |V_0(t)|,\ldots,|V_n(t)|,\ 
\theta_0(t),\ldots,\theta_n(t)\,\bigr]^\top .
\end{equation}
The state is equivalently represented by complex bus voltages. A fixed virtual-reference voltage sets the boundary magnitude and angle; a stiff branch anchors the solved physical substation voltage.

%=============================================================
\section{Proposed Methodology}
\label{sec:methodology}
\begin{figure}[t]
\centering
\includegraphics[width=0.95\linewidth]{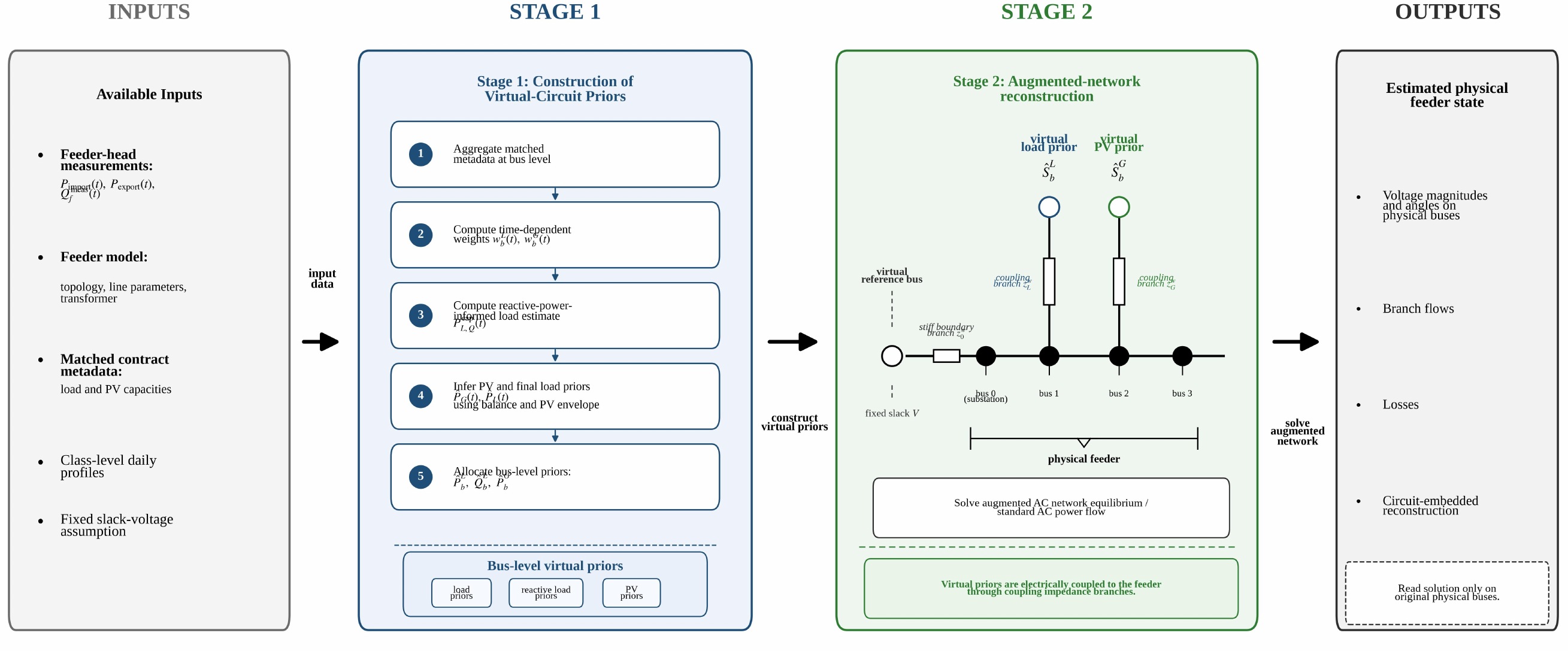}
\caption{Schematic of the proposed circuit-embedded reconstruction. Feeder-head measurements and matched metadata are first converted into bus-level virtual priors. The virtual prior buses are then connected to the physical feeder through coupling impedance branches, while the feeder boundary is represented by a stiff virtual slack connection.}
\label{fig:method_schematic}
\end{figure}

%=============================================================
The proposed reconstruction consists of two stages with different roles. Stage~1 uses the sparse information set to specify the active and reactive powers assigned to the virtual circuit elements. Stage~2 performs the state reconstruction by embedding those elements in the physical network and solving the augmented nonlinear AC equilibrium. The reconstructed state is its physical-bus component; Fig.~\ref{fig:method_schematic} summarizes the workflow.

\subsection{Stage 1: Construction of Virtual-Circuit Priors}
\label{subsec:stage1}
Let \(\mathcal B_L \subseteq \mathcal N\) denote the set of physical buses with matched consumption metadata, and let \(\mathcal B_G \subseteq \mathcal N\) denote the set of physical buses with matched downstream production metadata. Multiple contract records may be associated with the same physical bus. The implementation therefore aggregates metadata at the bus level. For each physical bus \(b\), define
\begin{equation}
C_b^L=\sum_{i\in\mathcal L_b} C_i^L,
\qquad
C_b^G=\sum_{g\in\mathcal G_b} C_g^G,
\end{equation}
where \(\mathcal L_b\) and \(\mathcal G_b\) denote the matched consumption and production records connected to bus \(b\), respectively. All capacity quantities used in the allocation formulas are expressed as active-power-equivalent proxies. The corresponding static capacity factors are
\begin{equation}
\alpha_b=\frac{C_b^L}{\sum_{\ell\in\mathcal B_L} C_\ell^L},
\qquad
\beta_b=\frac{C_b^G}{\sum_{h\in\mathcal B_G} C_h^G}.
\label{eq:allocation}
\end{equation}
These factors define static structural priors over the spatial distribution of downstream consumption and production. In the implementation, static factors are combined with predefined 24-point load and PV profiles. Let \(\pi_b^L(t)\) denote the assigned load profile normalized to unit daily mean, and let \(\pi_b^G(t)\) denote the assigned production profile divided by its daily peak, with values below \(0.10\) set to zero. The time-varying allocation weights are
\begin{equation}
w_b^L(t)
=
\frac{C_b^L\pi_b^L(t)}
{\sum_{\ell\in\mathcal B_L}C_\ell^L\pi_\ell^L(t)},
\qquad
w_b^G(t)
=
\frac{C_b^G\pi_b^G(t)}
{\sum_{h\in\mathcal B_G}C_h^G\pi_h^G(t)}.
\label{eq:time_weights}
\end{equation}
Because the load profiles and total matched load capacity are positive, the load-weight denominator is nonzero. If the PV-weight denominator vanishes, the PV envelope in~\eqref{eq:pv_prior_total} is zero, all production priors are set to zero, and the implementation uses uniform fallback PV weights, which
leave these zero-valued production priors unchanged. The aggregate active-load prior is constructed by combining a reactive-power-informed load estimate with an approximate active-power balance and a bounded PV envelope.
The reactive-power-informed estimate uses the measured feeder-head reactive power and an effective profile-weighted load power factor:
\begin{equation}
\bar{\tau}(t)
=
\sum_{b\in\mathcal B_L}
w_b^L(t)\tan\!\left(\arccos(\mathrm{pf}_b)\right),
\end{equation}
where \(\mathrm{pf}_b\) is the assumed class-level load power factor at bus \(b\). 

The reactive-power-informed active-load estimate is
\begin{equation}
P_{L,Q}(t)=
\left\{
\begin{array}{@{}l@{\;}l@{}}
Q_{\mathrm f}^{\mathrm{meas}}(t)/\bar{\tau}(t),
& Q_{\mathrm f}^{\mathrm{meas}}(t)\!>\!0
  \text{ and }\bar{\tau}(t)\!\ge\!\epsilon_\tau,\\
\max\{P_{\mathrm f}^{\mathrm{meas}}(t),0\},
& \text{otherwise},
\end{array}
\right.
\label{eq:pload_q}
\end{equation}
where \(\epsilon_\tau\) is a small numerical tolerance used to avoid excessive amplification when the effective reactive-to-active load ratio is close to zero. 
The reported estimator uses \(\mathrm{pf}=0.95\) for both residential and commercial loads, following a benchmark assumption used in studies of the IEEE European LV feeder~\cite{khan_reduced_2022}. Consequently, \(\bar\tau(t)=\tan(\arccos(0.95))\approx0.3287\) over the evaluated horizon, so the \(\epsilon_\tau=10^{-6}\) guard is inactive. To avoid unrealistic allocations when feeder-head reactive power would otherwise over-amplify the load estimate, define
\begin{equation}
P_{L,Q}^{\mathrm{cap}}(t)
=
\min\{P_{L,Q}(t),1.05\bar C_L\},
\qquad
\bar C_L=\sum_{\ell\in\mathcal B_L}C_\ell^L .
\label{eq:pload_q_cap}
\end{equation}
The aggregate PV prior is obtained from active-power balance subject to a profile-based production envelope:
\begin{equation}
\begin{aligned}
\hat P_G(t)
=
\min\{&
\max[(1+\eta)P_{L,Q}^{\mathrm{cap}}(t)
      -P_{\mathrm f}^{\mathrm{meas}}(t),0],\\
&\gamma\sum_{b\in\mathcal B_G}C_b^G\pi_b^G(t)
\},
\end{aligned}
\label{eq:pv_prior_total}
\end{equation}
where \(\eta\) is a fixed feeder-loss factor and \(\gamma\) is the aggregate PV-envelope factor. This expression follows the approximate active-power balance \(P_{\mathrm f}^{\mathrm{meas}}(t)\approx P_L(t)+P_{\mathrm{loss}}(t)-P_G(t)\), with \(P_{\mathrm{loss}}(t)\approx \eta P_L(t)\), under the convention that positive \(P_{\mathrm f}^{\mathrm{meas}}\) denotes net import at the feeder head. A second active-load prior is obtained by inverting the same approximate feeder-loss balance using the feeder-head active power and inferred production:
\begin{equation}
P_{L,P}(t)
=
\frac{\max\{P_{\mathrm f}^{\mathrm{meas}}(t)+\hat P_G(t),0\}}
{1+\eta}.
\label{eq:pload_p}
\end{equation}

The aggregate active-load prior is then
\begin{equation}
\hat P_L(t)
=
\min\left\{
\rho P_{L,Q}^{\mathrm{cap}}(t)
+
(1-\rho)P_{L,P}(t),
1.05\bar C_L
\right\},
\label{eq:pload_total}
\end{equation}
where \(\rho\in[0,1]\) weights the direct contribution of the reactive-power-informed active-load estimate. 
The reported configuration uses \(\rho=0\). When the PV prior in~\eqref{eq:pv_prior_total} lies strictly between zero and its envelope, substitution into~\eqref{eq:pload_p} yields \(P_{L,P}(t)=P_{L,Q}^{\mathrm{cap}}(t)\). Thus, in these snapshots, the reactive-power-informed estimate also determines the active-load prior through the approximate balance.
The two values can differ when a PV bound is active; the final capacity safeguard may further clip \(\hat P_L(t)\). Feeder-head reactive power therefore enters the active-power construction through the PV estimate in~\eqref{eq:pv_prior_total}, and it is also used directly in the reactive-power allocation below.
%No empirical active-load debiasing factor is applied.
The bus-level priors are finally assigned as
\begin{align}
\hat P_b^L(t)
&= w_b^L(t)\hat P_L(t),
&& b\in\mathcal B_L, \label{eq:bus_p_load}\\
\hat P_b^G(t)
&= w_b^G(t)\hat P_G(t),
&& b\in\mathcal B_G, \label{eq:bus_p_pv}\\
\hat Q_b^G(t)
&=0,
&& b\in\mathcal B_G. \label{eq:bus_q_pv}
\end{align}
The load reactive-power prior is assigned as
\begin{equation}
\hat Q_b^L(t)=
\begin{cases}
s_Q(t)\,\hat Q_{b,0}^L(t), & Q_{\mathrm f}^{\mathrm{meas}}(t)>0,\\
w_b^L(t)Q_{\mathrm f}^{\mathrm{meas}}(t), & Q_{\mathrm f}^{\mathrm{meas}}(t)\le 0,
\end{cases}
\label{eq:bus_q_load}
\end{equation}
where \(\hat Q_{b,0}^L(t)= \hat P_b^L(t)\tan(\arccos(\mathrm{pf}_b))\).
Here \(s_Q(t)\) is a bounded scaling factor used to align the aggregate inductive reactive-power prior with the measured feeder-head reactive power. For \(Q_{\mathrm f}^{\mathrm{meas}}(t)>0\), it is computed as
\begin{equation}
s_Q(t)=
\mathrm{clip}\!\left(
\frac{Q_{\mathrm f}^{\mathrm{meas}}(t)}
{\sum_{b\in\mathcal B_L}\hat Q_{b,0}^L(t)},
s_Q^{\min},s_Q^{\max}
\right),
\label{eq:q_scaling}
\end{equation}
where \(s_Q^{\min}\) and \(s_Q^{\max}\) are clipping bounds used to avoid excessive amplification of the reactive-power prior. The reported implementation uses \(s_Q^{\min}=0.80\), \(s_Q^{\max}=1.20\), \(\rho=0\), \(\eta=0.05\), and \(\gamma=0.70\). In the implementation, the scaling in~\eqref{eq:q_scaling} is applied only when \(\sum_{b\in\mathcal B_L}\hat Q_{b,0}^L(t)\) is strictly above a small numerical tolerance. Otherwise, signed feeder-head reactive power is allocated directly using the load weights, or a zero reactive prior is assigned when the measured reactive power is negligible.

\subsection{Stage 2: Augmented-Network Reconstruction}
\label{subsec:stage2}
Stage~2 embeds the priors in an augmented network and solves its nonlinear AC equilibrium; the reconstructed state is the physical-bus component. Smaller coupling impedances reduce virtual--physical voltage separation.
In the reported implementation, the downstream coupling resistances are \(2\times\) (load) and \(4\times\) (production) the median physical line-segment resistance, the coupling reactance is \(0.5\times\) the median line-segment reactance, and the boundary branch uses \(0.005\times\) the median resistance and reactance with a \(10^{-5}\,\Omega\) floor (Table~\ref{tab:params_results}). These impedances control virtual--physical voltage separation, branch losses, and current redistribution in the AC solution, while the assigned constant-power priors are held fixed. For each bus-level downstream prior, an auxiliary virtual bus is created. A virtual load-prior bus \(\tilde b_L\) is introduced for each \(b\in\mathcal B_L\), and a virtual production-prior bus \(\tilde b_G\) is introduced for each \(b\in\mathcal B_G\). Each virtual prior bus \(k\) is connected to its associated physical bus through a coupling branch with impedance \(z_k^{w}=r_k^{w}+jx_k^{w}\).
The substation coupling branch is chosen much stiffer than downstream prior branches. Let \(b(k)\) denote the physical bus associated with virtual prior bus \(k\). If \(\tilde V_k\) and \(V_{b(k)}\) denote the corresponding virtual and physical voltages, the coupling current is
\begin{equation}
I_k^{w}
=
\frac{\tilde V_k-V_{b(k)}}{z_k^{w}}.
\label{eq:Iweight}
\end{equation}
Thus, \(I_k^w\) denotes the branch current from the virtual prior bus toward the associated physical bus; the corresponding term in the physical-bus KCL below uses the opposite current direction.
We adopt the convention that load prior power is positive for withdrawal from the feeder, whereas production prior power is positive for injection into the feeder. Currents are written using the convention that positive \(I_k^{\mathrm{prior}}\) denotes a current withdrawal from the virtual bus. Thus, a load prior appears with positive sign, while a production prior appears with negative sign. 

For \(b\in\mathcal B_L\), define the complex load prior \(\hat S_b^L(t)=\hat P_b^L(t)+j\hat Q_b^L(t)\). The corresponding virtual-bus current term is
\begin{equation}
I_b^{L,\mathrm{prior}}
=
\frac{\bigl(\hat S_b^L(t)\bigr)^*}{\tilde V_{b,L}^*}.
\label{eq:Iprior_load}
\end{equation}

For \(b\in\mathcal B_G\), define \(\hat S_b^G(t)=\hat P_b^G(t)+j\hat Q_b^G(t)\), where \(\hat Q_b^G(t)=0\) is assumed in the absence of downstream inverter-reactive-power telemetry. The corresponding virtual-bus current term is
\begin{equation}
I_b^{G,\mathrm{prior}}
=
-\frac{\bigl(\hat S_b^G(t)\bigr)^*}{\tilde V_{b,G}^*},
\label{eq:Iprior_prod}
\end{equation}
where the minus sign reflects that production acts as an injection into the feeder rather than a withdrawal from it.
The augmented network equilibrium is defined by Kirchhoff's current law at both physical and virtual buses. All current-balance equations below are written with currents leaving the corresponding bus taken as positive. Let \(\mathcal K_b\) denote the set of downstream virtual prior buses attached to physical bus \(b\), and let \(\mathcal K=\bigcup_{b\in\mathcal N}\mathcal K_b\) denote the set of all downstream virtual prior buses. Let \(\mathcal N_b\) denote the set of physical buses adjacent to bus \(b\), and let \(Y_{bj}\) denote the corresponding branch admittance between buses \(b\) and \(j\). For each non-substation physical bus \(b\neq 0\), Kirchhoff's current law gives
\begin{equation}
\sum_{j\in\mathcal N_b} Y_{bj}\bigl(V_b - V_j\bigr)
+
\sum_{k\in\mathcal K_b}
\frac{V_b-\tilde V_k}{z_k^w}
=0.
\label{eq:physical_kcl}
\end{equation}

For the physical substation bus \(0\), the stiff boundary branch adds the term
\begin{equation}
\sum_{j\in\mathcal N_0} Y_{0j}\bigl(V_0 - V_j\bigr)
+
\sum_{k\in\mathcal K_0}
\frac{V_0-\tilde V_k}{z_k^w}
+
\frac{V_0-\tilde V_0}{z_0^w}
=0.
\label{eq:substation_kcl}
\end{equation}

Equations~\eqref{eq:physical_kcl}--\eqref{eq:substation_kcl} are written in compact branch form for clarity. In the implementation, the corresponding AC power-flow equations use the full \texttt{pandapower} network model, including transformer and admittance-matrix contributions. For compactness, \(I_k^{\mathrm{prior}}\) denotes either the load-prior current in~\eqref{eq:Iprior_load} or the production-prior current in~\eqref{eq:Iprior_prod}. At each downstream virtual prior bus \(k\in\mathcal K\), Kirchhoff's current law gives
\begin{equation}
\frac{\tilde V_k-V_{b(k)}}{z_k^w}
+
I_k^{\mathrm{prior}}
=0.
\label{eq:virtual_kcl}
\end{equation}
Equations~\eqref{eq:physical_kcl}--\eqref{eq:virtual_kcl}, together with the virtual-reference slack condition, define the Stage~2 AC equilibrium and enforce the prescribed virtual powers to numerical tolerance.

\subsubsection{Substation Boundary Anchoring}
A virtual reference bus with fixed voltage \(\tilde V_0=1.02\angle0^\circ\,\mathrm{p.u.}\), which provides the angle reference, is connected to physical bus \(0\) through \(z_0^w=r_0^w+jx_0^w\), with \(|z_0^w|\ll|z_k^w|\). Thus, the physical substation voltage is solved but stiffly anchored. The reference is fixed because time-varying boundary-voltage telemetry is unavailable, while feeder-head active and reactive power determine the aggregate operating level through Stage~1.

\subsection{Implementation}
In the augmented \texttt{pandapower} network~\cite{thurner_pandapoweropen-source_2018}, original downstream injections are disabled and bus-aggregated priors are reintroduced only at virtual buses. The external grid is placed at the virtual reference bus, the augmented network is solved by Newton--Raphson power flow, and the reconstructed state is read from the physical buses; coupling-branch flows are retained only as diagnostics. The implementation uses a balanced network equivalent, with one complex voltage per bus. The synthetic reference uses the same balanced representation, so the reported errors do not capture phase unbalance. The first snapshot uses a flat initialization and the remaining 671 use the preceding solution.
%=============================================================
\section{Case Study and Evaluation}
\label{sec:validation}
%=============================================================
\begin{figure}[t]
\centering
\includegraphics[width=0.90\linewidth]{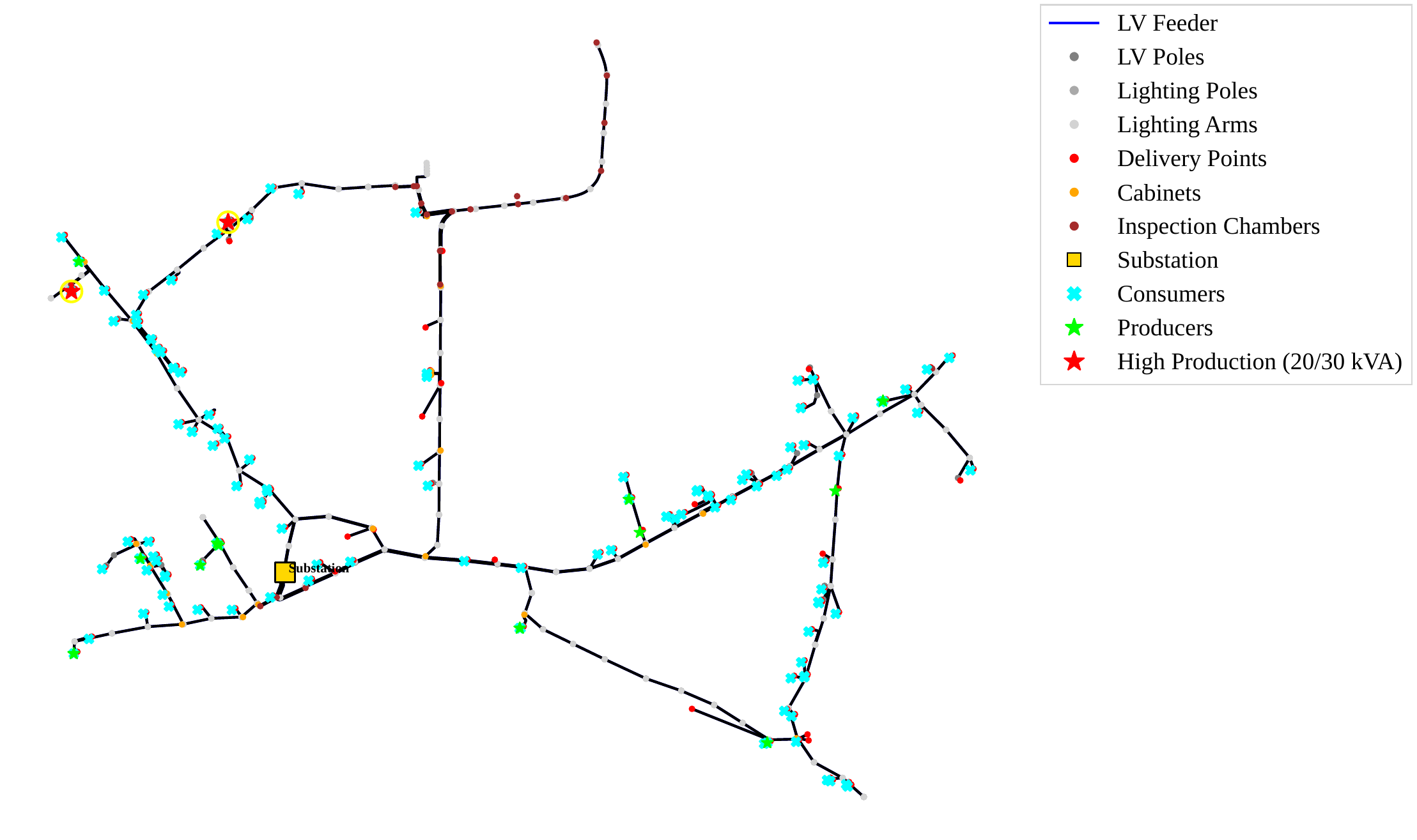}
\caption{Low-voltage feeder, showing the feeder infrastructure, matched consumer locations, matched production locations, and the feeder substation.}
\label{fig:topology}
\end{figure}
The proposed methodology is evaluated on a real Portuguese low-voltage distribution feeder in Vila Nova de Famalic\~ao, modeled in \texttt{pandapower} and observed under the targeted extreme-sparsity regime. The network, illustrated in Fig.~\ref{fig:topology}, consists of 459 buses and 601 lines, with a nominal low-voltage level of 230/400\,V. Physical feeder losses are computed by AC power flow using the feeder's cable impedances and segment lengths. Cable parameters are obtained from technical specifications, including~\cite{E-REDES_Derivações_2025}, and the resulting network is solved in \texttt{pandapower}~\cite{thurner_pandapoweropen-source_2018}.
The matched metadata contain 138 consumption-contract records assigned to 108 downstream load elements on 106 physical buses and 17 production-contract records assigned to 14 physical buses. Multiple colocated records are retained when forming capacity totals and then aggregated at the bus level, producing 106 load-prior branches and 14 production-prior branches in the augmented network. The matched capacity proxies total \(1.183350\,\mathrm{MW}\) for load and \(0.099230\,\mathrm{MW}\) for production. Metadata-driven profile classes use the 75th percentile of bus-level capacity, yielding load and PV thresholds of \(0.013800\,\mathrm{MW}\) and \(0.007208\,\mathrm{MW}\), respectively. Table~\ref{tab:params_results} summarizes the implementation parameters and reconstruction results.

\begin{table}[t]
\centering
\caption{Implementation parameters and synthetic-benchmark
results. E/T denotes estimated/true.}
\label{tab:params_results}
\setlength{\tabcolsep}{1pt}
\renewcommand{\arraystretch}{0.90}
\scriptsize
\resizebox{\columnwidth}{!}{%
\begin{tabular}{@{}lc@{\hspace{3pt}}lc@{}}
\toprule
\multicolumn{2}{c}{\textbf{Implementation}} &
\multicolumn{2}{c}{\textbf{Reconstruction}} \\
\cmidrule(r){1-2}\cmidrule(l){3-4}
\textbf{Parameter} & \textbf{Value} &
\textbf{Metric} & \textbf{Value} \\
\midrule
Slack voltage & \(1.0200\) p.u. &
Solved snapshots & \(672/672\) \\

\(\eta/\gamma/\rho\) & \(0.05/0.70/0.00\) &
Valid entries & \(308{,}448/308{,}448\) \\

\(s_Q^{\min}/s_Q^{\max}\) & \(0.80/1.20\) &
Voltage RMSE & \(0.0033\) p.u. \\

\(\epsilon_\tau\) & \(10^{-6}\) &
Voltage MAE & \(0.0014\) p.u. \\

Load power factor & \(0.95\) &
Max \(|e_V|\) & \(0.0337\) p.u. \\

Median \(r/x\) & \(0.010219/0.000869\,\Omega\) &
Mean \(r_P\) & \(0.614\) kW \\

Load branch \(r\) & \(0.020438\,\Omega\) &
Fraction \(r_P>0\) & \(0.3423\) \\

PV branch \(r\) & \(0.040876\,\Omega\) &
Mean \(V_{\min}\), E/T & \(0.9253/0.9306\) p.u. \\

Coupling \(x\) & \(0.000435\,\Omega\) &
Min. \(V_{\min}\), E/T & \(0.7834/0.8156\) p.u. \\

Boundary \(r/x\) & \(0.000051/0.000010\,\Omega\) &
Mean summed \(|\)transfer\(|\) & \(167.985\) kW \\
\bottomrule
\end{tabular}%
}
\end{table}
\subsection{Evaluation Methodology}
Because downstream ground-truth telemetry is unavailable, evaluation uses a calibrated synthetic benchmark of 672 consecutive 15-minute snapshots from October~1--7, 2023. Truth injections used separate peak-normalized residential and commercial weekday/weekend profiles. Commercial classes were assigned at or above the 75th percentile of element-level capacity proxies, including zero-capacity elements in the load threshold. Unmatched load elements had zero capacity proxies and therefore zero active and reactive truth powers. Capacities and time-invariant asset scales were independently perturbed by \(\mathcal U[0.95,1.05]\) factors. Load and PV utilization factors were \(0.22\) and \(0.70\), load power factor was \(0.95\), and PV output was capped at \(95\%\) of perturbed capacity. Common daily load and PV multipliers were drawn from \(\mathcal U[0.95,1.05]\) and \(\mathcal U[0.80,1.00]\), respectively. Independent asset-level snapshot multipliers followed \(\mathcal N(1,0.01^2)\) for loads and \(\mathcal N(1,0.015^2)\) for PV, clipped to \([0.80,1.25]\) and \([0.70,1.20]\), respectively. No feeder-head measurement noise or PV time shift was introduced. Each truth snapshot was obtained by AC power flow and then reduced to the information available to the reconstruction method: feeder-head active/reactive power, fixed boundary voltage, network model, and matched metadata. The estimator uses a distinct but structurally related profile library and does not observe realized nodal profiles, perturbations, or internal voltages. Truth and estimator share the nominal load power factor, zero PV reactive power, and \(1.02\)-p.u. boundary voltage, while truth PV utilization and the estimator envelope factor \(\gamma\) both equal \(0.70\). Moreover, \(\eta,\gamma,\rho\) were assessed on the same time window. The reported errors are therefore in-sample and exclude boundary-voltage error. All snapshots solved without back-off or exclusion; errors are pooled over the \(308{,}448\) physical-bus/time voltage-magnitude pairs.
Across the 672 snapshots, the aggregate Stage~1 load prior had MAE, RMSE, and signed bias of \(3.047\), \(3.704\), and \(+1.422\,\mathrm{kW}\), respectively. The corresponding PV-prior values were \(1.894\), \(2.929\), and \(+1.481\,\mathrm{kW}\). Mean true aggregate load and PV production were \(146.404\) and \(18.678\,\mathrm{kW}\), respectively; the prior RMSEs correspond to \(2.53\%\) and \(15.68\%\) of these means. These aggregate errors were measured under the shared load power factor and zero PV reactive power, conditions that favor the reactive-power-informed estimate \(Q_{\mathrm f}^{\mathrm{meas}}/\bar\tau\). With \(\rho=0\), substitution gives \(\hat P_L(t)=P_{L,Q}^{\mathrm{cap}}(t)\) in the 163 snapshots with an interior PV prior. 

\subsection{Feeder-Head Power Consistency}
\begin{figure}[t]
\centering
\begin{minipage}[t]{0.48\linewidth}
\centering
\includegraphics[width=\linewidth]{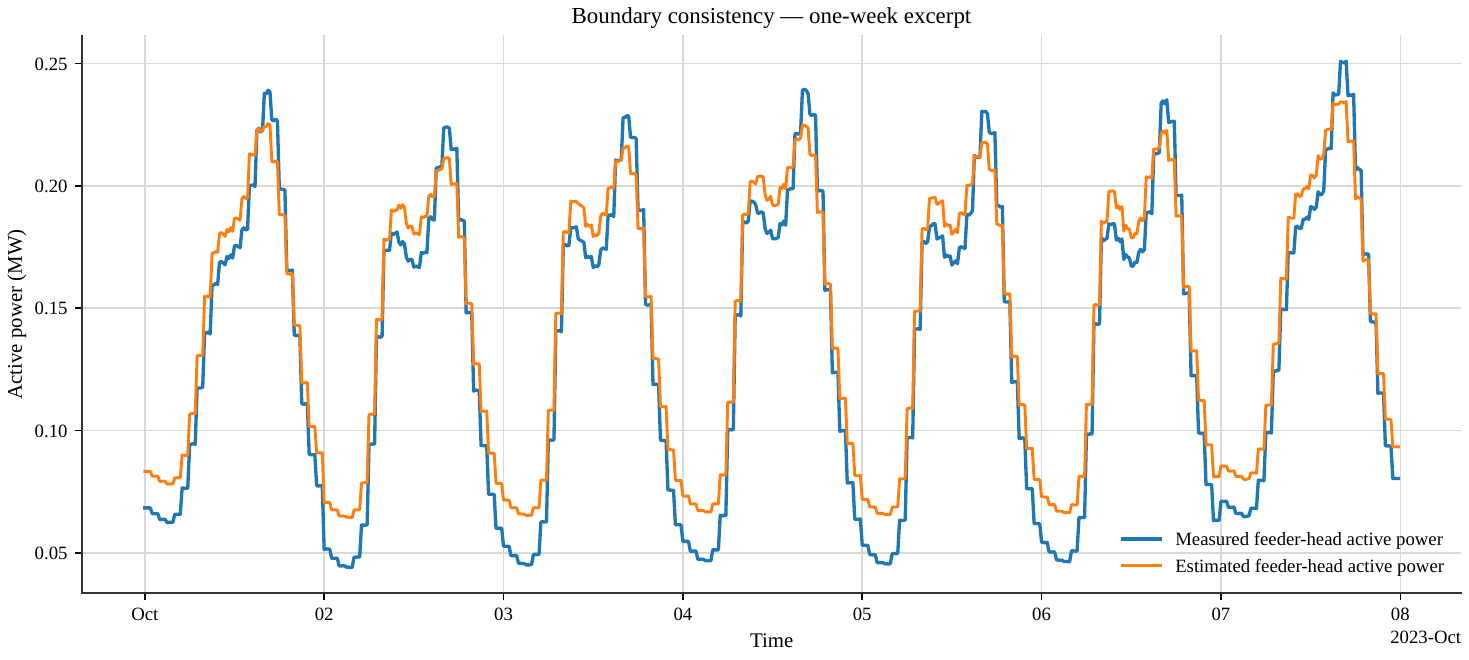}
\captionof{figure}{Feeder-head active-power consistency over the one-week benchmark.}
\label{fig:boundary}
\end{minipage}
\hfill
\begin{minipage}[t]{0.48\linewidth}
\centering
\includegraphics[width=\linewidth]{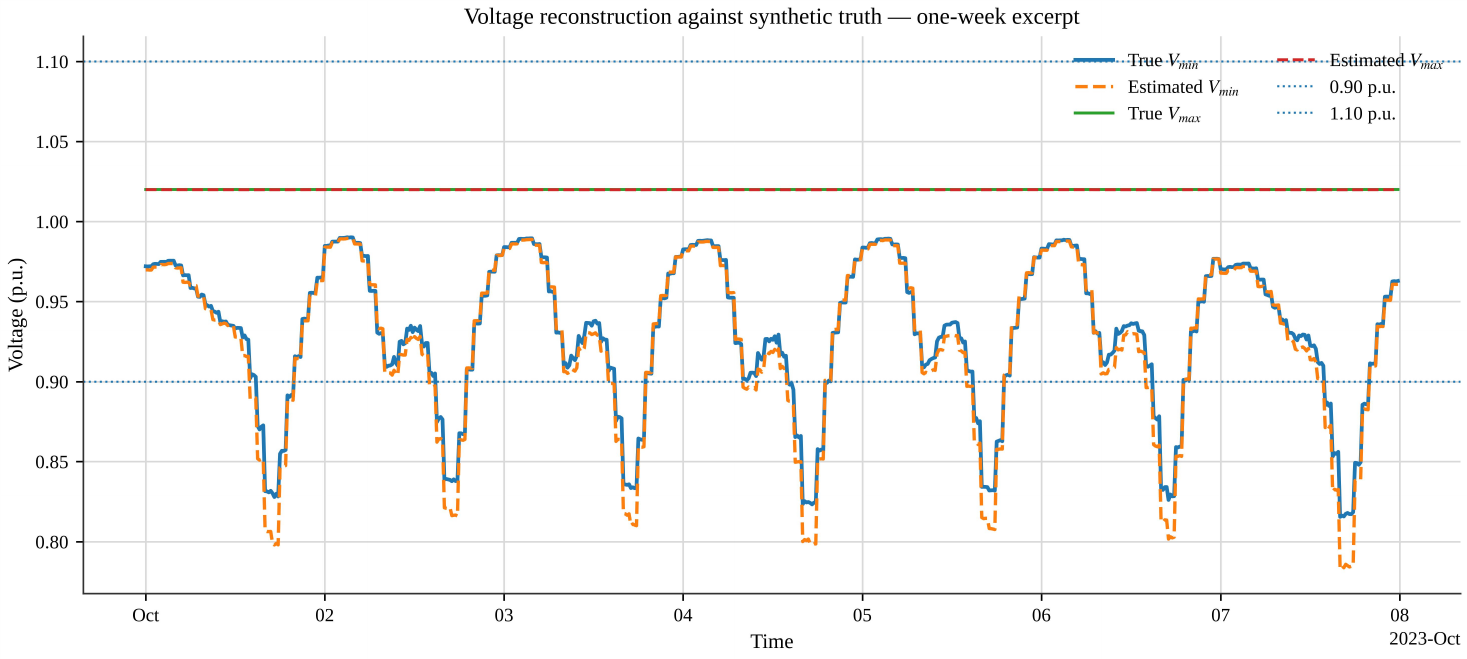}
\captionof{figure}{Voltage reconstruction against synthetic truth.}
\label{fig:true_est_voltage}
\end{minipage}
\end{figure}

The signed feeder-head active-power residual is defined as
\begin{equation}
r_P(t)=P_{\mathrm f}^{\mathrm{rec}}(t)-P_{\mathrm f}^{\mathrm{meas}}(t),
\label{eq:feeder_residual}
\end{equation}
where \(P_{\mathrm f}^{\mathrm{rec}}(t)\) is the active power supplied at the virtual-slack terminal of the boundary branch, including its losses. The mean signed feeder-head active-power residual was \(0.614\,\mathrm{kW}\), or \(0.45\%\) of the mean true feeder-head active power. The maximum absolute residual was \(15.394\,\mathrm{kW}\), and \(34.2\%\) of snapshots had a positive signed residual. The corresponding reactive-power residuals are also small, with mean signed residual \(1.107\,\mathrm{kVAr}\) and maximum absolute residual \(3.379\,\mathrm{kVAr}\). No reactive-scaling clamp is activated. %; since \(\sum_{b\in\mathcal B_L}\hat Q_{b,0}^L(t)=\bar\tau(t)\hat P_L(t)\), \(s_Q(t)=1\) whenever \(\hat P_L(t)=P_{L,Q}^{\mathrm{cap}}(t)=Q_{\mathrm f}^{\mathrm{meas}}(t)/\bar\tau(t)\). 
Fig.~\ref{fig:boundary} is a feeder-head consistency diagnostic. Feeder-head active power informs the priors but is not imposed as an equality constraint in Stage~2, so nonzero boundary residuals can occur.
\subsection{Voltage Reconstruction and Numerical Convergence}
Fig.~\ref{fig:true_est_voltage} shows that the reconstructed voltage trajectory remained close to the synthetic truth. The solver converged for all 672 snapshots without rejection, with mean reconstructed \(V_{\min}=0.9253\,\mathrm{p.u.}\) versus \(0.9306\,\mathrm{p.u.}\) in the reference trajectory; \(V_{\min}<0.90\,\mathrm{p.u.}\) occurred in 130 truth snapshots and 156 reconstructed snapshots. Of these, 128 coincided; the reconstruction reported 28 snapshots below \(0.90\,\mathrm{p.u.}\) that were not below this level in the truth, and missed two truth snapshots. 

%\subsection{Coupling-Branch Diagnostics}
%At the virtual-bus terminals, the temporal mean of the summed absolute active-power deviations from the assigned priors was \(1.14\times10^{-5}\,\mathrm{kW}\); the mean ratio of this deviation to total absolute transferred active power was \(8.12\times10^{-8}\). The mean and maximum sums of absolute transferred active power were \(167.985\) and \(315.589\,\mathrm{kW}\). 

\subsection{Direct-Prior AC Power-Flow Consistency Check}
As an implementation check, we applied the same Stage~1 priors directly at physical buses, without virtual prior buses or coupling branches. Both cases solve AC circuit equilibria. Direct-prior and augmented-circuit voltage RMSEs were \(0.003189\) and \(0.003293\) p.u., respectively; their voltage states differed by \(1.74\times10^{-4}\) p.u. RMSE and \(1.06\times10^{-3}\) p.u. at most. Since \(\tilde V_k-V_{b(k)}=z_k^w I_k^w\), the embedding approaches direct injection as \(|z_k^w|\to0\) for bounded coupling currents. The selected impedances caused a small voltage-state perturbation and slightly increased RMSE relative to direct-prior AC power flow. At the virtual-bus terminals, the mean summed absolute deviation from the assigned priors was \(1.14\times10^{-5}\,\mathrm{kW}\).

\subsection{Comparison with WLS Baselines}
\begin{figure}[t]
\centering
\begin{minipage}[t]{0.48\linewidth}
\vspace{0pt}
\centering
\captionof{table}{Circuit-equilibrium reconstruction and WLS baselines.}
\label{tab:circuit_wls_comparison}
\vspace{1pt}
\setlength{\tabcolsep}{0.8pt}
\renewcommand{\arraystretch}{0.90}
\tiny
\begin{tabular}{@{}lcccc@{}}
\toprule
\textbf{Metric} & \textbf{Circuit} & \textbf{All V}
& \textbf{Slack V} & \textbf{Matched} \\
\midrule
Conv. (\%)        & 100.00 & 100.00 & 100.00 & 100.00 \\
RMSE              & 0.0033 & 0.0457 & 0.0156 & 0.0029 \\
MAE               & 0.0014 & 0.0317 & 0.0068 & 0.0015 \\
Max \(|e_V|\)     & 0.0337 & 0.2115 & 0.1146 & 0.0262 \\
Mean \(V_{\min}\) & 0.9253 & 1.0029 & 0.8979 & 0.9261 \\
\bottomrule
\end{tabular}
\end{minipage}
\hfill
\begin{minipage}[t]{0.48\linewidth}
\vspace{8pt}
\centering
\includegraphics[width=\linewidth]{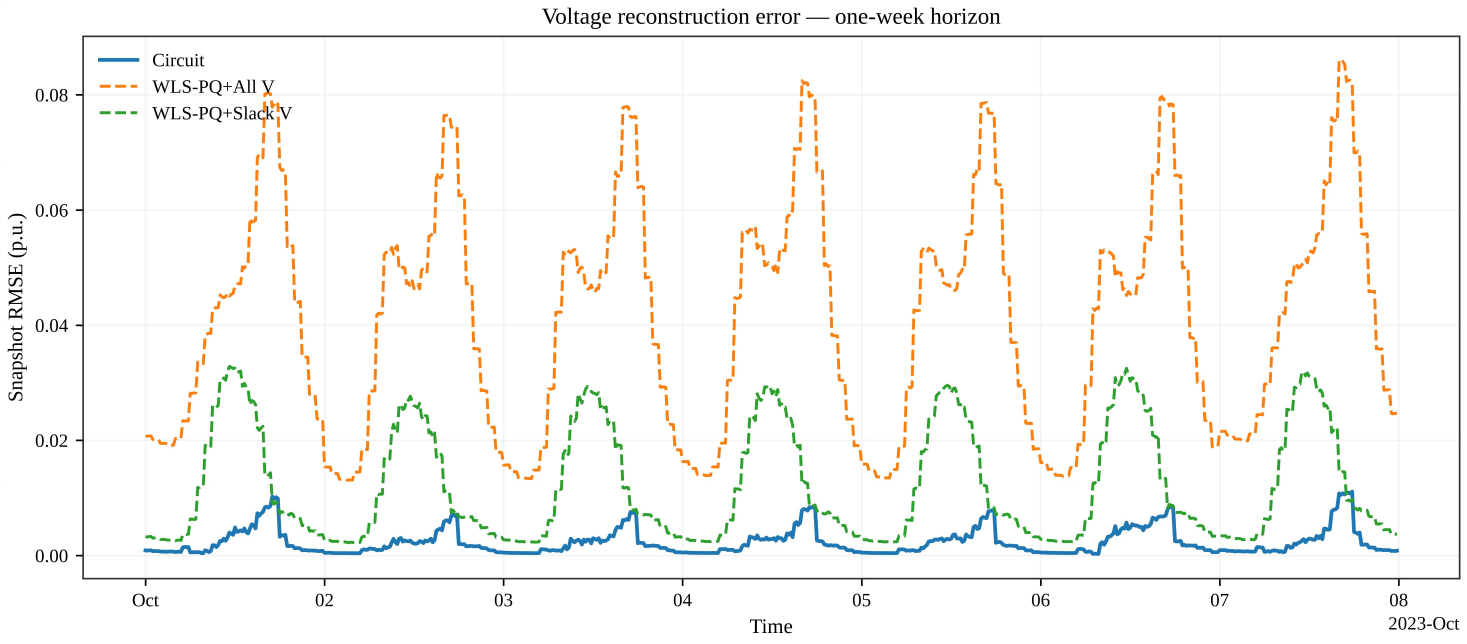}
\captionof{figure}{Voltage RMSE for the circuit reconstruction and the
two static-prior WLS baselines; matched-prior WLS is reported in Table~\ref{tab:circuit_wls_comparison}.}
\label{fig:wls_circuit_rmse}
\end{minipage}
\end{figure}

Two static-prior WLS baselines use the same topology, feeder-head measurements, fixed slack voltage, and matched metadata. Both include the slack-voltage measurement and bus-level active/reactive pseudo-measurements; WLS-PQ+All V additionally includes \(1.02\)-p.u. voltage pseudo-measurements at non-slack buses with standard deviation \(0.05\,\mathrm{p.u.}\) Their nodal pseudo-measurements use static allocation factors, as in conventional sparse-observability practice, and omit the profile information used in Stage~1. Bus \(P/Q\) pseudo-measurements use \(50\%\) relative standard deviations, with respective lower bounds of \(0.010\,\mathrm{MW}\) and \(0.010\,\mathrm{MVAr}\).
The 327 structurally passive buses are enforced as zero-injection constraints in all WLS variants. Zero-valued pseudo-measurements at buses containing injection equipment retain standard deviations of \(0.050\,\mathrm{MW}\) and \(0.050\,\mathrm{MVAr}\).
Feeder-head \(P/Q\) measurements use \(1\%\) relative standard deviations with respective lower bounds of \(0.005\,\mathrm{MW}\) and \(0.005\,\mathrm{MVAr}\), while the slack-voltage standard deviation is \(10^{-4}\,\mathrm{p.u.}\)
For the nodal injection pseudo-measurements, define \(a_b^L=\alpha_b\) for \(b\in\mathcal B_L\) and \(a_b^L=0\) otherwise; similarly, define \(a_b^G=\beta_b\) for \(b\in\mathcal B_G\) and \(a_b^G=0\) otherwise. With positive bus \(P/Q\) denoting net withdrawal, the static-prior WLS pseudo-measurements are
\begin{align}
z_{P,b}^{\mathrm{WLS}}(t)
&=
P_{\mathrm{import}}(t)a_b^L
-
P_{\mathrm{export}}(t)a_b^G,\\
z_{Q,b}^{\mathrm{WLS}}(t)
&=
Q_{\mathrm f}^{\mathrm{meas}}(t)a_b^L .
\end{align}
When \(P_{\mathrm{export}}(t)=0\), the static-prior WLS pseudo-measurements assign no separate downstream production component. 
Extending load priors by zero outside \(\mathcal B_L\) and production priors by zero outside \(\mathcal B_G\), a separate matched-prior WLS mechanism check uses the Stage~1 outputs:
\begin{align}
z_{P,b}^{\mathrm{prof}}(t)
&=
\hat P_b^L(t)-\hat P_b^G(t),\\
z_{Q,b}^{\mathrm{prof}}(t)
&=
\hat Q_b^L(t).
\end{align}
Nodal \(P/Q\) pseudo-measurements are assigned only at non-slack buses. The slack bus carries the feeder-head measurements \(z_{P,0}=-P_{\mathrm f}^{\mathrm{meas}}\) and \(z_{Q,0}=-Q_{\mathrm f}^{\mathrm{meas}}\), consistent with the net-withdrawal sign convention.
Matched-prior WLS uses the same Stage~1 outputs as the augmented circuit, providing a comparison under common priors. Table~\ref{tab:circuit_wls_comparison} and Fig.~\ref{fig:wls_circuit_rmse} summarize the comparison; all methods converged for all snapshots. 
Circuit reconstruction achieved RMSE \(0.003293\,\mathrm{p.u.}\), versus \(0.015642\) for WLS-PQ+Slack V and \(0.045675\,\mathrm{p.u.}\) for WLS-PQ+All V, reductions of \(78.9\%\) and \(92.8\%\), respectively. WLS-PQ+All V had mean \(V_{\min}=1.0029\,\mathrm{p.u.}\), versus \(0.9306\,\mathrm{p.u.}\) in the truth. Matched-prior WLS achieved RMSE \(0.002886\,\mathrm{p.u.}\), MAE \(0.001516\,\mathrm{p.u.}\), and maximum error \(0.026230\,\mathrm{p.u.}\) It achieved lower RMSE and maximum error than circuit reconstruction, while circuit reconstruction retained lower MAE (\(0.001392\) versus \(0.001516\,\mathrm{p.u.}\)).
The static-prior comparisons evaluate the complete procedure, while the common-prior comparisons assess the reconstruction formulation and measurement reconciliation. 
Direct calibration against downstream measurements is unavailable in this field regime. The circuit formulation requires no pseudo-measurement variances, while WLS standard deviations are computed at each snapshot using the prescribed relative-error, minimum-value, and zero-value rules.
%Because downstream measurements are unavailable in field operation, pseudo-measurement variances cannot be calibrated against nodal data; the circuit formulation requires none, whereas the WLS results above rely on variances fixed a priori.

\subsection{Sensitivity to Calibration Parameters}
All tested configurations converged. Relative to the base RMSE of \(0.003293\,\mathrm{p.u.}\), setting \(\eta=0\) or \(\gamma=0.40\) increased RMSE to \(0.008531\) and \(0.008120\,\mathrm{p.u.}\), respectively. The truth-implied mean loss factor was approximately \(0.050\), while the truth PV-utilization factor and selected \(\gamma\) were both \(0.70\); these results therefore constitute in-sample parameter assessment. Increasing \(\rho\) from \(0\) to \(1\) raised RMSE from \(0.003293\) to \(0.003686\,\mathrm{p.u.}\) and the mean signed active-power residual from \(0.614\) to \(2.610\,\mathrm{kW}\), supporting the selected \(\rho=0\). The PV envelope was zero in 336 snapshots and its upper bound was active in 173; the zero-PV lower bound was never active with a positive envelope. These 509 snapshots were exactly those in which changing \(\rho\) altered the Stage~1 aggregate load prior. The prior was unchanged in the remaining 163 snapshots.
%=============================================================
\section{Conclusion}
\label{sec:conclusion}
%=============================================================
Under boundary-only telemetry, the method constructs bus-level priors from feeder-head powers, contract metadata, and profiles, then embeds them as constant-power elements in the AC circuit equations. In the in-sample synthetic benchmark with shared load power factor and zero PV reactive power, all 672 snapshots converged; voltage RMSE, MAE, and maximum error were \(0.0033\), \(0.0014\), and \(0.0337\) p.u., respectively. The complete procedure outperformed two static-prior WLS baselines, while matched-prior WLS achieved lower RMSE \(0.0029\) versus \(0.0033\,\mathrm{p.u.}\) Direct-prior AC power flow differed from the augmented circuit by only \(1.74\times10^{-4}\,\mathrm{p.u.}\) RMSE. The selected coupling impedances therefore produced a small voltage-state perturbation. The common-prior comparisons indicate that the improvement over the selected static-prior WLS baselines primarily comes from Stage~1 prior construction. The AC embedding requires no pseudo-measurement variances, avoiding their specification in the reconstruction stage. 

\section*{Acknowledgment}
LLM tools (ChatGPT~\cite{openai_chatgpt}, Claude~\cite{anthropic_claude}, and Gemini~\cite{google_gemini}) assisted with review, proofreading, and consistency checks. The authors retain full responsibility for all content and results.
%=============================================================
% References
%=============================================================

\bibliographystyle{ieeetr}
\bibliography{references}
\end{document}